\documentclass[12pt,twoside]{article}

\usepackage[latin1]{inputenc}
\usepackage{graphicx}
\usepackage{color}
\usepackage[longnamesfirst,round]{natbib}

\usepackage{enumitem}
\usepackage{latexsym}
\usepackage{amssymb}
\usepackage{epsfig}
\usepackage{hyperref}
\usepackage{xcolor}

\usepackage[paper=a4paper, top=2cm, bottom=3cm, left=2.5cm, right=3.0cm]{geometry} 

\usepackage{amsmath}
\usepackage{float}
\usepackage{dsfont}
\usepackage{multirow}
\usepackage{bm}

\allowdisplaybreaks 

\def\bX{\mathbf{X}}
\def\bx{\mathbf{x}}

\newtheorem{theo}{Theorem}[section]
\newtheorem{lemma}[theo]{Lemma}
\newtheorem{rem}[theo]{Remark}

\begin{document}

\title{\bf Identification in a Fully Nonparametric Transformation Model with Heteroscedasticity
}

\author{{\sc\Large Nick Kloodt}\footnote{Bundesstra\ss{}e 55, 20146 Hamburg, Nick.Kloodt@uni-hamburg.de, 040 428387167, ORCID 0000-0001-8998-1703}\\ Department of Mathematics, University of Hamburg}

\maketitle

\begin{abstract}\noindent

The so far most general identification result in the context of nonparametric transformation models is proven. The result is constructive in the sense that it provides an explicit expression of the transformation function.

\end{abstract}


\section{Introduction} 
\def\theequation{1.\arabic{equation}}
\setcounter{equation}{0}

The underlying question of this article can be formulated quite easily: Given some real valued random variable $Y$ and some $\mathbb{R}^{d_X}$-valued random variable $\bX$ fulfilling the heteroscedastic transformation model
\begin{equation}\label{trafoeq}
h(Y)=g(\bX)+\sigma(\bX)\varepsilon
\end{equation}
with some error term $\varepsilon$ independent of $\bX$ and fulfilling $E[\varepsilon]=0$ and $\operatorname{Var}(\varepsilon)=1$, are the model components $h:\mathbb{R}\rightarrow\mathbb{R},\ g:\mathbb{R}^{d_X}\rightarrow\mathbb{R},\sigma:\mathbb{R}^{d_X}\rightarrow(0,\infty)$ and the error distribution uniquely determined if the joint distribution of $(Y,\bX)$ is known? This uniqueness is called identification of a model.

Over the last years, transformation models have attracted more and more attention since they are often used to obtain desirable properties by first transforming the dependent random variable of a regression model. Applications for such transformations can reach from reducing skewness of the data to inducing additivity, homoscedasticity or even normality of the error terms. Already \citet{BC1964}, \citet{BD1981} and \citet{ZR1969} introduced some parametric classes of transformation functions. \citet{Hor1996} proved for a linear regression function $g$ and homoscedastic errors that the model is identified, when $h(y_0)=0$ is assumed for some $y_0\in\mathbb{R}$ and the regression parameter is standardized such that the first component, which is different from zero, is equal to one. Later, the ideas of \citet{Hor1996} were extended by \citet{EHN2004} to general smooth regression functions $g$. The arguably most general identification results so far were provided by \citet{CKC2015} and \citet{VvK2018}, who considered general regression functions and homoscedastic errors as well, but allowed endogenous regressors. \citet{LSvK2008} used similar ideas to obtain identifiability of a model with parametric transformation functions as a special case. As will be seen in Section \ref{idea} these approaches can not be applied to the heteroscedastic model so that different methods are needed. Despite their practical relevance (e.g. in duration models, see \citet{KST2011}), results allowing heteroscedasticity are rare. \citet{ZLJ2009} showed identifiability in a single-index model with a linear regression function $g$ and a known variance function $\sigma^2$. \citet{WW2018} applied this model to lung cancer data. \citet{NNvK2016} required identifiability implicitly in their assumptions.

In contrast to the approaches mentioned above, it is tried here to avoid any parametric assumption on $h,g$ or $\sigma$, which to the author's knowledge has not been done before. Note that the validity of the model is unaffected by linear transformations. This means that for arbitrary constants $a>0,b\in\mathbb{R}$ equation (\ref{trafoeq}) still holds when replacing $h$, $g$ and $\sigma$ by
$$\tilde{h}(y)=ah(y)+b,\qquad\tilde{g}(\bx)=ag(\bx)+b\qquad\textup{and}\qquad\tilde{\sigma}(\bx)=a\sigma(\bx).$$
Of course, one could have chosen an arbitrary $a\in\mathbb{R}$ as well, but similar to existing results the transformation function $h$ will be restricted to be strictly increasing without loss of generality. Nevertheless, at least two conditions for fixing $a$ and $b$ are needed. Referring to the fact that these conditions will determine the linear transformation they are sometimes called location and scale constraints.

This remainder is organized as follows. First, some assumptions are listed before the main identification result for heteroscedastic transformation models is motivated and stated. Afterwards, a short conclusion in Section \ref{conclusion} is followed by the Appendix, which contains some results on uniqueness of solutions to differential equations and the proof of the main result.

\section{The Idea and the Result}\label{idea}

Before the identification result can be motivated, some assumptions and notations have to be introduced. First, basic assumptions concerning validity of model (\ref{trafoeq}) and continuity of its model components are given.
\begin{enumerate}[label=(\textbf{A\arabic{*}})]
	\item\label{A1} Let $Y,\varepsilon$ and $\bX$ be real valued and $\mathbb{R}^{d_X}$-valued random variables, respectively, with
					$$h(Y)=g(\bX)+\sigma(\bX)\varepsilon$$
					for some transformation, regression and variance functions $h,g$ and $\sigma^2$.
	\item\label{A2} $\varepsilon$ is a centred random variable independent of $\bX$ with $E[\varepsilon]=0$ and $\operatorname{Var}(\varepsilon)=1$.
	\item\label{A3} Let the density $f_{\varepsilon}$ of $\varepsilon$ be continuous and let $h,g$ and $\sigma$ from \ref{A1} be continuously differentiable.
\end{enumerate}
Moreover, a regularity assumption for the conditional distribution function $F_{Y|\bX}$ of $Y$ given $\bX$ is needed.
\begin{enumerate}[label=(\textbf{A4})]
	\item\label{A4} The conditional cumulative distribution function $(y,\bx)\mapsto F_{Y|\bX}(y|\bx)$ is continuously differentiable with respect to $y$ and $\bx$. Let $v\geq0$ be a weight function with support $\operatorname{supp}(v)$ such that $\frac{\partial}{\partial y}F_{Y|\bX}(y|\bx)>0$ for all $y\in\mathbb{R},\bx\in\operatorname{supp}(v)$ and such that (with $g$ and $\sigma$ from \ref{A1})
	$$A:=\int v(\bx)\left(\frac{\sigma(\bx)\frac{\partial g(\bx)}{\partial x_1}-g(\bx)\frac{\partial\sigma(\bx)}{\partial x_1}}{\sigma(\bx)}\right)\,d\bx\quad\textup{and}\quad B:=\int v(\bx)\frac{\frac{\partial\sigma(\bx)}{\partial x_1}}{\sigma(\bx)}\,d\bx$$
	are well defined with $B\neq0$.
\end{enumerate}
The assumption $B\neq0$ requires heteroscedasticity of the model. Note that the homoscedastic case was already treated by \citet{CKC2015}. Later, it will be shown in Remark \ref{remtesthomoscedasticity} that \ref{A1}--\ref{A3} and the first part of \ref{A4} exclude the case that there exist a homoscedastic and a heteroscedastic version of model (\ref{trafoeq}) at the same time. In the following, the functions $h,g,\sigma$ and $f_{\varepsilon}$ from \ref{A1} and \ref{A3} are used to show their uniqueness and consequently identification of the model.

\noindent
\subsection{The Transformation Function as a Solution to an Initial Value Problem}
Many of the homoscedastic identification approaches mentioned in the introduction are based on the same idea (see \citet{EHN2004}, \citet{Hor2009} and recently \citet{CKC2015}). Using the example of \citet{CKC2015} their method can be summarized in the following way: Let $F_{Y|\bX}$ be the conditional cumulative distribution function of $Y$ conditioned on $\bX$. Take the derivatives of $F_{Y|\bX}$ with respect to $y$ and some component of $\bx$, divide the first by the latter one and obtain the transformation function by integrating this quotient. After applying some identification constraints the transformation function is identified as it only depends on the joint distribution of $(Y,\bX)$. In heteroscedastic models, the reasoning has to be changed since the way, the transformation function enters the conditional distribution function and its partial derivatives, becomes more complex. The latter functions can be written as
\begin{align}
F_{Y|\bX}(y|\bx)&=P(Y\leq y|\bX=\bx)\nonumber
\\[0,2cm]&=P\bigg(\varepsilon\leq\frac{h(Y)-g(\bX)}{\sigma(\bX)}\Big|\bX=\bx\bigg)\nonumber
\\[0,2cm]&=F_{\varepsilon}\bigg(\frac{h(y)-g(\bx)}{\sigma(\bx)}\bigg),\nonumber
\\[0,5cm]\frac{\partial F_{Y|\bX}(y|\bx)}{\partial y}&=f_{\varepsilon}\bigg(\frac{h(y)-g(\bx)}{\sigma(\bx)}\bigg)\frac{h'(y)}{\sigma(\bx)}>0\label{idFYX}
\end{align}
and
$$\frac{\partial F_{Y|\bX}(y|\bx)}{\partial x_i}=-f_{\varepsilon}\bigg(\frac{h(y)-g(\bx)}{\sigma(\bx)}\bigg)\frac{\sigma(\bx)\frac{\partial g(\bx)}{\partial x_i}+(h(y)-g(\bx))\frac{\partial\sigma(\bx)}{\partial x_i}}{\sigma(\bx)^2},\quad i=1,...,d_X.$$
Here, $h'(y)$ is an abbreviation for the derivative $\frac{\partial}{\partial y}h(y)$ and $F_{\varepsilon}$ denotes the cumulative distribution function of $\varepsilon$. Hence, even if \ref{A4} is valid the transformation function can not be obtained by simply integrating the quotient
\begin{equation}\label{lnxinverse}
\frac{\frac{\partial F_{Y|\bX}(y|\bx)}{\partial y}}{\frac{\partial F_{Y|\bX}(y|\bx)}{\partial x_i}}=-\frac{h'(y)\sigma(\bx)}{\sigma(\bx)\frac{\partial g(\bx)}{\partial \bx_i}+(h(y)-g(\bx))\frac{\partial\sigma(\bx)}{\partial x_i}},
\end{equation}
since the denominator now also depends on the transformation function.

Instead, we consider the reciprocal value of (\ref{lnxinverse}) and name this $\tilde{\lambda}$:
$$\tilde{\lambda}(y|\bx):=\frac{\frac{\partial F_{Y|\bX}(y|\bx)}{\partial x_i}}{\frac{\partial F_{Y|\bX}(y|\bx)}{\partial y}}=-\frac{\sigma(\bx)\frac{\partial g(\bx)}{\partial \bx_i}+(h(y)-g(\bx))\frac{\partial\sigma(\bx)}{\partial x_i}}{h'(y)\sigma(\bx)}.$$
Next, if $v$ is the weight function from \ref{A4} $\tilde{\lambda}$ can be integrated with respect to $\bx$ as follows to obtain
\begin{equation}\label{lambda}
\lambda(y):=\int v(\bx)\lambda(y|\bx)\,d\bx\,=-\frac{A+Bh(y)}{h'(y)}
\end{equation}
with $A$ and $B\neq0$ from \ref{A4}. Since assumption \ref{A4} implies $h'>0$ and consequently strict monotonicity of $h$, there exists exactly one root of $\lambda$ which will be called
$$y_0:=\lambda^{-1}(0)$$
in the following. Due to (\ref{lambda}) it holds that $h(y_0)=-\frac{A}{B}$.

In the following, the problem of identifying model (\ref{trafoeq}) is reduced to solving an ordinary differential equation uniquely. Afterwards, basic uniqueness theorems for initial value problems will imply the main identification result. To this end, rewrite equation (\ref{lambda}) to obtain
\begin{equation}\label{diffeq}
h'(y)=-\frac{A+Bh(y)}{\lambda(y)}
\end{equation}
for all $y\in(y_0,\infty)$. This indeed can be understood as a differential equation, but an initial condition is needed to obtain an initial value problem. Here, the initial condition
\begin{equation}\label{scaleconstraint}
h(y_1)=\alpha
\end{equation}
for some $y_1>y_0$ and some $\alpha>-\frac{A}{B}$ is considered (remember that $h$ was assumed to be strictly increasing). Theorem \ref{diffuni} in the appendix yields uniqueness of any solution to this initial value problem on any interval $[z_1,z_2]\subseteq(y_0,\infty)$. This identification result can be generalized to all $y\in\mathbb{R}$.

\noindent
\subsection{Uniqueness of the Unknown Coefficients}\label{uniquecoefficients}
The reasoning above is designed for fixed $A$ and $B$, that is, it remains to prove uniqueness of these coefficients. Moreover, it would be desirable to derive an explicit formula for the transformation function instead of only proving its uniqueness. This will be done in the remainder of this section.

First, the initial value problem, which corresponds to the equations (\ref{diffeq}) and (\ref{scaleconstraint}), is solved by
\begin{equation}\label{solutionAB}
h(y)=\frac{(A+B\alpha)\exp\bigg(-B\int_{y_1}^y\frac{1}{\lambda(u)}\,du\bigg)-A}{B}\quad\textup{for all }y\in(y_0,\infty).
\end{equation}
By straightforward calculations, it can be verified that (\ref{solutionAB}) is indeed a solution to the initial value problem.
Second, as was already mentioned in the introduction, model (\ref{trafoeq}) is not only fulfilled for $h,g,\sigma$, but also for any linear transformation of these functions. Therefore, to obtain uniqueness it is necessary to fix these linear transforms. This can be done by requiring so called location and scale constraints and corresponds to fixing $A$ and $\alpha$. While
\begin{equation}\label{locconstraint}
h(y_0)=0
\end{equation}
is chosen as the location constraint the scale constraint is equal to the initial condition (\ref{scaleconstraint}), that is, $\alpha$ is viewed as an arbitrary, but fixed positive number. Here, the location  constraint was chosen such that equation (\ref{lambda}) implies $A=0$. Nevertheless, other location constraints are conceivable as well as can be seen in Remark \ref{otherlocconstraints}.

Consequently, equation (\ref{solutionAB}) reduces to
\begin{equation}\label{solutionB}
h(y)=\alpha\exp\bigg(-B\int_{y_1}^y\frac{1}{\lambda(u)}\,du\bigg)\quad\textup{for all }y\in(y_0,\infty).
\end{equation}
If there exist two coefficients $B\neq\tilde{B}$ such that the corresponding transformation functions from (\ref{solutionB}) fulfil model (\ref{trafoeq}), it would hold that
$$\tilde{h}(y)=\alpha \bigg(\frac{h(y)}{\alpha}\bigg)^{\frac{\tilde{B}}{B}}\quad\textup{for all }y\in(y_0,\infty).$$
Assume without loss of generality $\tilde{B}>B>0$. Then,
$$\tilde{h}'(y)=\frac{\tilde{B}}{B}\bigg(\frac{h(y)}{\alpha}\bigg)^{\frac{\tilde{B}}{B}-1}h'(y)\overset{y\searrow y_0}{\longrightarrow}0.$$
Therefore, continuous differentiability of $h$ and $\tilde{h}$ would imply $\tilde{h}'(y_0)=0$, which due to (\ref{idFYX}) would lead to a violation of \ref{A4}. Hence, $B$ is unique under \ref{A1}--\ref{A4}, which finally leads to the main identification result. Note that the same argument is valid for transformation functions as in (\ref{solutionAB}) since these are simply linearly transformed versions of (\ref{solutionB}).
\begin{theo}\label{theoid}
	Let $y_2<y_0<y_1$ and assume \ref{A1}--\ref{A4} and (\ref{scaleconstraint}).
	\begin{enumerate}[label=\alph{*})]
		\item For each $A\in\mathbb{R}$ such that $\alpha>-\frac{A}{B}$, the unique solution to (\ref{diffeq}) on $(y_0,\infty)$ is given by (\ref{solutionAB}). It can be extended to a global unique solution to (\ref{diffeq}) by\vspace{0,2cm}
		\begin{equation}\label{globalsolution}
		h(y)=\left\{\begin{array}{ll}
		\frac{(A+B\alpha)\exp\Big(-B\int_{y_1}^y\frac{1}{\lambda(u)}\,du\Big)-A}{B}&y>y_0\\
		-\frac{A}{B}&y=y_0\\
		\frac{(A+B\alpha_2)\exp\Big(-B\int_{y_2}^y\frac{1}{\lambda(u)}\,du\Big)-A}{B}&y<y_0
		\end{array}\right.,\vspace{0,3cm}
		\end{equation}
		where $\alpha_2$ is uniquely determined by requiring $\underset{y\searrow y_0}{\lim}\,h'(y)=\underset{y\nearrow y_0}{\lim}\,h'(y)=h'(y_0)$ as
		\begin{equation}\label{deflambdaunid}
		\alpha_2=-\frac{\underset{t\rightarrow0}{\lim}\,(A+B\alpha)\exp\bigg(B\bigg(\int_{y_2}^{y_0-t}\frac{1}{\lambda(u)}\,du-\int_{y_1}^{y_0+t}\frac{1}{\lambda(u)}\,du\bigg)\bigg)+A}{B}.
		\end{equation}
		\item If additionally (\ref{locconstraint}) and $\alpha=1$ hold, one has
		\begin{equation}\label{globalsolutionidconstraints}
		h(y)=\left\{\begin{array}{ll}
		\exp\Big(-B\int_{y_1}^y\frac{1}{\lambda(u)}\,du\Big)&y>y_0\\
		0&y=y_0\\
		\alpha_2\exp\Big(-B\int_{y_2}^y\frac{1}{\lambda(u)}\,du\Big)&y<y_0\end{array}\right.,\vspace{0,3cm}
		\end{equation}
		where $\alpha_2$ is uniquely determined by requiring $\underset{y\searrow y_0}{\lim}\,h'(y)=\underset{y\nearrow y_0}{\lim}\,h'(y)=h'(y_0)$ as
		\begin{equation}\label{deflambda2}
		\alpha_2=-\underset{t\rightarrow0}{\lim}\,\exp\bigg(B\bigg(\int_{y_2}^{y_0-t}\frac{1}{\lambda(u)}\,du-\int_{y_1}^{y_0+t}\frac{1}{\lambda(u)}\,du\bigg)\bigg).
		\end{equation}
	\end{enumerate}
	Moreover, $B$ is uniquely determined and it holds that
	$$g(\bx)=E[h(Y)|\bX=\bx]\quad\textup{and}\quad\sigma(\bx)=\sqrt{\operatorname{Var}(h(Y)|\bX=\bx)}.$$
\end{theo}\noindent
\noindent
The proof can be found in Section \ref{proofextensionh}.

Finally, two remarks are given dealing on the one hand with further generalizations and implications for future estimation and testing techniques and on the other hand with justifying alternative identification constraints.
\begin{rem}\label{remtesthomoscedasticity}
	\begin{enumerate}[label=\alph{*})]
		\item If $\sigma$ is not constant, there are values $\bx\in\mathbb{R}^{d_X}$ such that $\frac{\partial}{\partial x_1}\sigma(\bx)\neq0$. Consequently, $\tilde{\lambda}$ changes its sign for these $\bx$. If $\sigma$ is constant, that is, the error is homoscedastic, this is not the case. Hence, model (\ref{trafoeq}) can not be fulfilled for homoscedastic and heteroscedastic errors at the same time.
		\item The identification result can be generalized in many regards. For example, one could have used any other partial derivative $\frac{\partial}{\partial x_i},i\neq1,$ in \ref{A4} as well. Moreover, $v$ can be chosen as a Dirac delta function as well and it is possible to consider error densities $f_{\varepsilon}$ with bounded support. See \citet{Klo2019} for a more detailed examination. Moreover, it is conjectured that the result can be generalized to conditional independence of $\bX$ and $\varepsilon$ given endogenous regressors similarly to \citet{CKC2015}.
	\end{enumerate}
\end{rem}
\begin{rem}\label{otherlocconstraints}
	One could have used other scale and location constraints than (\ref{scaleconstraint}) and (\ref{locconstraint}). For example, consider for some real numbers $y_a<y_b,\alpha_a<\alpha_b$ the conditions
	\begin{equation}\label{otherconstraints}
	h(y_a)=\alpha_a\quad\textup{and}\quad h(y_b)=\alpha_b.
	\end{equation}
	Assume there exist two transformation functions $T_1,T_2$ such that model (\ref{trafoeq}) and the constraints (\ref{otherconstraints}) are fulfilled. Then, the functions
	$$\tilde{T}_1(y)=\frac{T_1(y)-T_1(y_0)}{T_1(y_1)-T_1(y_0)}\quad\textup{and}\quad\tilde{T}_2(y)=\frac{T_2(y)-T_2(y_0)}{T_2(y_1)-T_2(y_0)}$$
	fulfil model (\ref{trafoeq}) and the constraints (\ref{scaleconstraint}) and (\ref{locconstraint}). This leads to $\tilde{T}_1\equiv\tilde{T}_2$ so that
	$$T_1(y)=(\alpha_b-\alpha_a)\frac{\tilde{T}_1(y)-\tilde{T}_1(y_a)}{\tilde{T}_1(y_b)-\tilde{T}_1(y_a)}+\alpha_a=(\alpha_b-\alpha_a)\frac{\tilde{T}_2(y)-\tilde{T}_2(y_a)}{\tilde{T}_2(y_b)-\tilde{T}_2(y_a)}+\alpha_a=T_2(y)$$
	for all $y\in\mathbb{R}$.
	
	A similar reasoning can be applied to show that identification constraints like
	$$h(y_a)=\alpha_a,\quad h'(y_a)=\alpha_b$$
	for some $y_a,\alpha_a\in\mathbb{R},\alpha_b>0$ ensure uniqueness of the transformation function as well.
\end{rem}
\noindent

\noindent
\section{Conclusion and Outlook}\label{conclusion}
The so far most general identification result in the theory of transformation models has been provided. While doing so, the techniques of \citet{EHN2004} and \citet{CKC2015} have been used to reduce the problem of identifiability to that of solving an ordinary differential equation. Most of the previous results are contained as special cases. The main contribution consists in allowing heteroscedastic errors, which justifies the common practice to assume identifiability like for example in the paper of \citet{NNvK2016}.

Moreover, the result is constructive in the sense that it does not only guarantee identification of the model, but even supplies an analytic expression of the transformation function depending on the joint cumulative distribution function of the data and some parameter $B$. This parameter is identified, too, and can be expressed as in \citet{Klo2019} under the additional assumption of a twice continuously differentiable transformation function.

Due to the explicit character of equation (\ref{solutionAB}), future research could consist in analysing the resulting plug-in estimator. This will be the topic of a subsequent paper. Furthermore, the presented results could be successively generalized as in Remark \ref{remtesthomoscedasticity} or by allowing vanishing derivatives of $h$. Moreover, it would be desirable to develop conditions on the joint distribution function of $(Y,X)$ under which model (\ref{trafoeq}) is fulfilled. In contrast to the thoughts on identifiability here, such a question addresses the solvability of (\ref{trafoeq}), that is, the issue of existence of a solution instead of uniqueness.

\appendix

\section{Uniqueness of Solutions to Ordinary Differential Equations}
In this Section, two basic results about ordinary differential equations and uniqueness of possible solutions are given. Theorem \ref{diffuni} is slightly modified compared to the version of \citet[p.~102]{For1999} so that the proof is presented as well.
\begin{lemma}\label{gronwall}(Gronwall's Inequality, see \citet{Gro1919} or \citet{Bel1953} for details)
	Let $I=[a,b]\subseteq\mathbb{R}$ be a compact interval. Let $u,v:I\rightarrow\mathbb{R}$ and $q:I\rightarrow[0,\infty)$ be continuous functions. Further, let
	$$u(y)\leq v(y)+\int_a^yq(z)u(z)\,dz$$
	for all $y\in I$. Then, one has
	$$u(y)\leq v(y)+\int_a^yv(z)q(z)\exp\bigg(\int_z^yq(t)\,dt\bigg)\,dz\quad\textup{for all }y\in I.$$
\end{lemma}\noindent

\begin{theo}\label{diffuni}(see \citet[p.~102]{For1999} for a related version)
	Let $b>a>y_0$ and $G\subseteq(y_0,\infty)\times\mathbb{R}^+$ be a set such that $[a,b]\times\mathbb{R}^+\subseteq G$. Moreover, let $D:G\rightarrow\mathbb{R},(y,h)\mapsto D(y,h),$ be continuous with respect to both components and continuously differentiable with respect to the second component. Then, for all $\theta_0>0$ any solution $h\in\mathcal{C}([a,b],\mathbb{R}^+)$ of the initial value problem
	$$h'(y)=D(y,h(y)),\quad h(a)=\theta_0$$
	is unique.
\end{theo}\noindent
\textbf{Proof:} Let $h_1,h_2:[a,b]\rightarrow\mathbb{R}^+$ be two solutions of the mentioned initial value problem. Since
$$K:=\{(y,\theta)\in[a,b]\times\mathbb{R}^+:y\in[a,b],\theta\in\{h_1(y),h_2(y)\}\}$$
is compact, there exists some $L>0$ such that $|D(y,\theta)-D(y,\psi)|\leq L|\theta-\psi|$ for all $(y,\theta),(y,\psi)\in K$. Consider the distance $d(y):=|h_1(y)-h_2(y)|$. Then for all $y\in[a,b]$
\begin{align*}
d(y)&=|h_1(y)-h_1(a)-(h_2(y)-h_2(a))|
\\[0,2cm]&=\bigg|\int_a^y(D(z,h_1(z))-D(z,h_2(z)))\,dz\bigg|
\\[0,2cm]&\leq\int_a^y|D(z,h_1(z))-D(z,h_2(z))|\,dz
\\[0,2cm]&\leq L\int_a^y|h_1(z)-h_2(z)|\,dz
\\[0,2cm]&=L\int_a^yd(z)\,dz.
\end{align*}
Gronwall's Inequality leads to $d\leq0$ (set $u=d,v\equiv0,q\equiv L$).\hfill$\square$

\section{Proof of Theorem \ref{theoid}}\label{proofextensionh}
Consider a compact interval $\mathcal{K}=[k_1,k_2]\subseteq(y_0,\infty)$ and recall equation (\ref{solutionAB}). Assumption \ref{A4} ensures $h'>0$. First, it is shown that $h$ as defined in (\ref{solutionAB}) is the unique solution to (\ref{diffeq}) on $[k_1,k_2]$. For the moment assume $k_1=y_1$ and define
$$G=[y_1,k_2]\times[\alpha,\infty)\quad\textup{and}\quad D:G\rightarrow\mathbb{R},\ D(y,h)=-\frac{A+Bh}{\lambda(y)}.$$
With the choices $a=y_1,b=k_2$ and $\theta_0=\alpha$, Theorem \ref{diffuni} ensures uniqueness of the solution to
$$h'(y)=-\frac{A+Bh(y)}{\lambda(y)}\quad\textup{for all }y\in[y_1,k_2].$$
By straightforward calculations, it can be verified that (\ref{solutionAB}) is indeed a solution to this initial value problem. Since $f_{\varepsilon}(y)>0$ for all $y\in\mathbb{R}$, this solution holds for arbitrarily large $k_2$. Hence, by letting $k_2$ tend to infinity uniqueness of $h$ on $[y_1,\infty)$ is obtained.

Now, consider an arbitrary value $\tilde{y}\in(y_0,y_1)$. Then, if the previous initial condition is replaced by
$$h(\tilde{y})=\tilde{\alpha}$$
for some $\tilde{\alpha}>0$, the same reasoning as before can be used to show that the differential equation
$$h'(y)=-\frac{A+Bh(y)}{\lambda(y)}\quad\textup{for all }y\in[\tilde{y},k_2]$$
is uniquely solved under this constraint by
\begin{align*}
h(y)&=\frac{(A+B\tilde{\alpha})\exp\Big(-B\int_{\tilde{y}}^y\frac{1}{\lambda(u)}\,du\Big)-A}{B}
\\[0,2cm]&=\frac{\frac{(A+B\tilde{\alpha})}{(A+B\alpha)}\exp\Big(-B\int_{\tilde{y}}^{y_1}\frac{1}{\lambda(u)}\,du\Big)(A+Bh(y))-A}{B}
\end{align*}
for all $y\in[\tilde{y},\infty)$, where the last equation follows from (\ref{solutionAB}). To fulfil the previous scale constraint $h(y_1)=\alpha$ it is required that
$$\tilde{\alpha}=\frac{(A+\alpha B)\exp\bigg(B\int_{\tilde{y}}^{y_1}\frac{1}{\lambda(u)}\,du\bigg)-A}{B}.$$
Since this in turn results in expression (\ref{solutionAB}) for all $y\in[\tilde{y},\infty)$, $h$ is identified for all $y\in[\tilde{y},\infty)$. Choosing $\tilde{y}$ arbitrarily close to $y_0$ results in
$$h(y)=\frac{(A+B\alpha)\exp\Big(-B\int_{y_1}^y\frac{1}{\lambda(u)}\,du\Big)-A}{B}\quad\textup{for all }y>y_0.$$
When proceeding analogously for $y<y_0$ with the initial condition
$$h(y_2)=\alpha'$$
for some $\alpha'<-\frac{A}{B}$, one has
$$h(y)=\frac{(A+B\alpha')\exp\Big(-B\int_{y_2}^y\frac{1}{\lambda(u)}\,du\Big)-A}{B}\quad\textup{for all }y<y_0.$$
Recall $h'(y)>0$ for all $y\in\mathbb{R}$ and let $t>0$. Due to the continuous differentiability of $h$ in $y_0$, one has
$$\underset{t\rightarrow0}{\lim}\,\frac{\frac{\displaystyle h(y_0+t)-h(y_0)}{\displaystyle t}}{\frac{\displaystyle h(y_0-t)-h(y_0)}{\displaystyle -t}}\rightarrow1.$$
On the other hand, it holds that
\begin{align*}
\frac{\frac{\displaystyle h(y_0+t)-h(y_0)}{\displaystyle t}}{\frac{\displaystyle h(y_0-t)-h(y_0)}{\displaystyle -t}}&=-\frac{(A+B\alpha)\exp\Big(-B\int_{y_1}^{y_0+t}\frac{1}{\lambda(u)}\,du\Big)}{(A+B\alpha')\exp\Big(-B\int_{y_2}^{y_0-t}\frac{1}{\lambda(u)}\,du\Big)}
\\[0,2cm]&=-\frac{(A+B\alpha)}{(A+B\alpha')}\exp\bigg(B\bigg(\int_{y_2}^{y_0-t}\frac{1}{\lambda(u)}\,du-\int_{y_1}^{y_0+t}\frac{1}{\lambda(u)}\,du\bigg)\bigg),
\end{align*}
so that
$$\alpha'=-\frac{\underset{t\rightarrow0}{\lim}\,(A+B\alpha)\exp\bigg(B\bigg(\int_{y_2}^{y_0-t}\frac{1}{\lambda(u)}\,du-\int_{y_1}^{y_0+t}\frac{1}{\lambda(u)}\,du\bigg)\bigg)+A}{B}=\alpha_2.$$
This leads to the uniqueness of solution (\ref{globalsolution}), since uniqueness of $B$ was already shown in Section \ref{uniquecoefficients}.\\
Inserting $A=0$ yields the second part of the assertion, while identification of $g$ and $\sigma$ as the conditional mean and standard deviation follows from standard arguments.\hfill$\square$

\section*{Acknowledgements}
This work was supported by the DFG (Research Unit FOR 1735 {\it Structural Inference in Statistics: Adaptation and Effciency}).\\
Moreover, I would like to thank Natalie Neumeyer and Ingrid Van Keilegom for their very helpful suggestions and comments on the project.

\bibliographystyle{plainnat}
\bibliography{References}

\end{document}